\documentclass[a4paper,12pt,leqno]{article}

\newcommand{\zero}{\setcounter{equation}{0}}

    \newcommand{\Om}{\Omega}
    \newcommand{\p}{\partial}
    \newcommand{\be}{\begin{equation}}
    \newcommand{\ee}{\end{equation}}
    \newcommand{\bea}{\begin{eqnarray}}
    \newcommand{\eea}{\end{eqnarray}}
    \newcommand{\bean}{\begin{eqnarray*}}
    \newcommand{\eean}{\end{eqnarray*}}
    \newcommand{\ba}{\begin{array}{ll}}
    \newcommand{\ea}{\end{array}}

\usepackage{amssymb,amsmath,amsthm,graphicx}

\newtheorem{thm}{Theorem}
\newtheorem{lem}{Lemma}

\begin{document}
 \title{A classical approach to dynamics of parabolic competitive systems}
\author{Bogdan Przeradzki\footnote{Address: Institute of Mathematics,
Technical University of \L \'od\'z,
ul. W\'olcza\'nska 215,
93-005  \L \'od\'z, Poland
}}

 \maketitle

\begin{abstract}
We study the reaction-diffusion system, its stationary solutions,
the behavior of the system near them and discuss similarities and
differences for different boundary conditions.
\end{abstract}



\section{Introduction}

We deal with a nonlinear parabolic system of the form
\begin{equation}
u_t=D\Delta u+f(u),    \label{sys1}
\end{equation}
with Neumann homogeneous boundary condition
\begin{equation}
\frac{\p u}{\p \nu}\left | _{\p\Om}=0   \right.    \label{neum}
\end{equation}
or, incidentally, Dirichlet homogeneous boundary condition
\begin{equation}
u\left | _{\p\Om}=0,  \right.       \label{dir}
\end{equation}
in a bounded domain $\Om\subset\Bbb{R}^m$ with  boundary  $\p\Om$
being $m-1$-dimensional sufficiently smooth manifold. Here, we
consider only classical solutions $u:\overline{\Om}\to \Bbb{R}^n $
(we emphasize that we have vector-valued functions, since
(\ref{sys1}) is, in fact, a system of $n$ parabolic equations. We
denote  by $D:={\rm diag}(d_1,\ldots , d_n)$ a diagonal matrix with
positive diagonal entries $d_i,$ $i=1,\ldots ,n,$  and by
$f:\Bbb{R}^n_+\to \Bbb{R}^n$ a continuous function defined on the
cone of nonnegative vectors $x=(x_1,\ldots ,x_n)\in\Bbb{R}^n,$ where
$x_i\ge 0,$ $i=1,\ldots ,n.$  Such parabolic systems are largely
investigated since they model kinetics of chemical reactions -- each
coordinate of $u$ measures density one of interacting components, or
they model development of several biological species living on the
same area $\Om $ and interacting  in different ways (preys and
predators, symbiosis, competition). We keep in mind the second
application that causes Neumann condition is more natural (comp.
\cite{Ch}). System (\ref{sys1}-\ref{neum}) defines a semiflow $\Phi
_t,$ $t\ge 0,$ on an appropriate space $X_{\alpha},$ if one uses the
theory of sectorial operators (comp. \cite{H,I}) or a semiflow  on a
cone of nonnegative continuous functions, if one applies the theory
of monotone dynamical systems (comp. \cite{S1}). We work in spaces
of continuous functions, since our main example of $f=(f_1,\ldots
,f_n)$ is
\begin{equation}
f_i(u)=u_i\left(1-\sum _{j=1}^n a_{ij}u_j\right),     \label{nonl}
\end{equation}
where all coefficients $a_{ij}$ are positive. If $f$ is a $C^1$-function such that
$$f_i(u)=0\quad {\rm for} \quad u_i=0,\quad i=1,\ldots ,n,$$
and
\begin{equation}
\frac{\p f_i}{\p u_j}\le 0, \quad {\rm for} \quad  i\neq j,   \label{compet}
\end{equation}
then (\ref{sys1}) defines a semiflow $\Phi _t,$ $t\ge 0,$ on the cone of nonnegative continuous functions $u:\overline{\Om}
\to [0,\infty)^n$ which is competitive in the sense of Hirsch (see \cite{S1,S2}). In particular, it means that all solutions of the system starting
with nonnegative functions $\varphi =u(\cdot ,0)\ge 0 $ are global in time and nonnegative for any $t>0.$ Moreover,  if
$$ u(\cdot ,T)\le \bar{u}(\cdot ,T)$$
are two such solutions comparable at any time $T>0,$ then
$$ u(\cdot ,t)\le \bar{u}(\cdot ,t)\quad {\rm for} \quad t<T.$$

\section{Steady-state solutions, single species}
\zero

First, we are interested in steady-state solutions, i.e. time
independent solutions. They satisfy elliptic system:
\begin{equation}
 \Delta u=-D^{-1}f(u),\quad   \frac{\p u}{\p \nu}\left | _{\p\Om}=0   \right. .  \label{ell}
\end{equation}
The obvious examples of such functions are zeros of the nonlinear term $f.$  An important question is if there are not exist
other steady-state solutions. They correspond nonuniform distributions of populations on the environment $\Om $ which do not change
in time.
We can exactly investigate the case of one species $(n=1).$ If $\Om\subset \Bbb{R},$ the analysis is standard and rather simple.
\begin{equation}
u''=-D^{-1}u(1-u) \quad {\rm in}\quad (0,L),\quad u'(0)=0=u'(L),\quad u\ge 0.     \label{steady}
\end{equation}
The second order ODE is an example of a conservative system with one degree of freedom so it has a first integral (energy)
$$E(u,u')=\frac{u'^2}{2}+\frac{u^2}{2D}-\frac{u^3}{3D}.$$
It enables to find the phase portrait of the system:

\includegraphics[angle=-90,width=12cm]{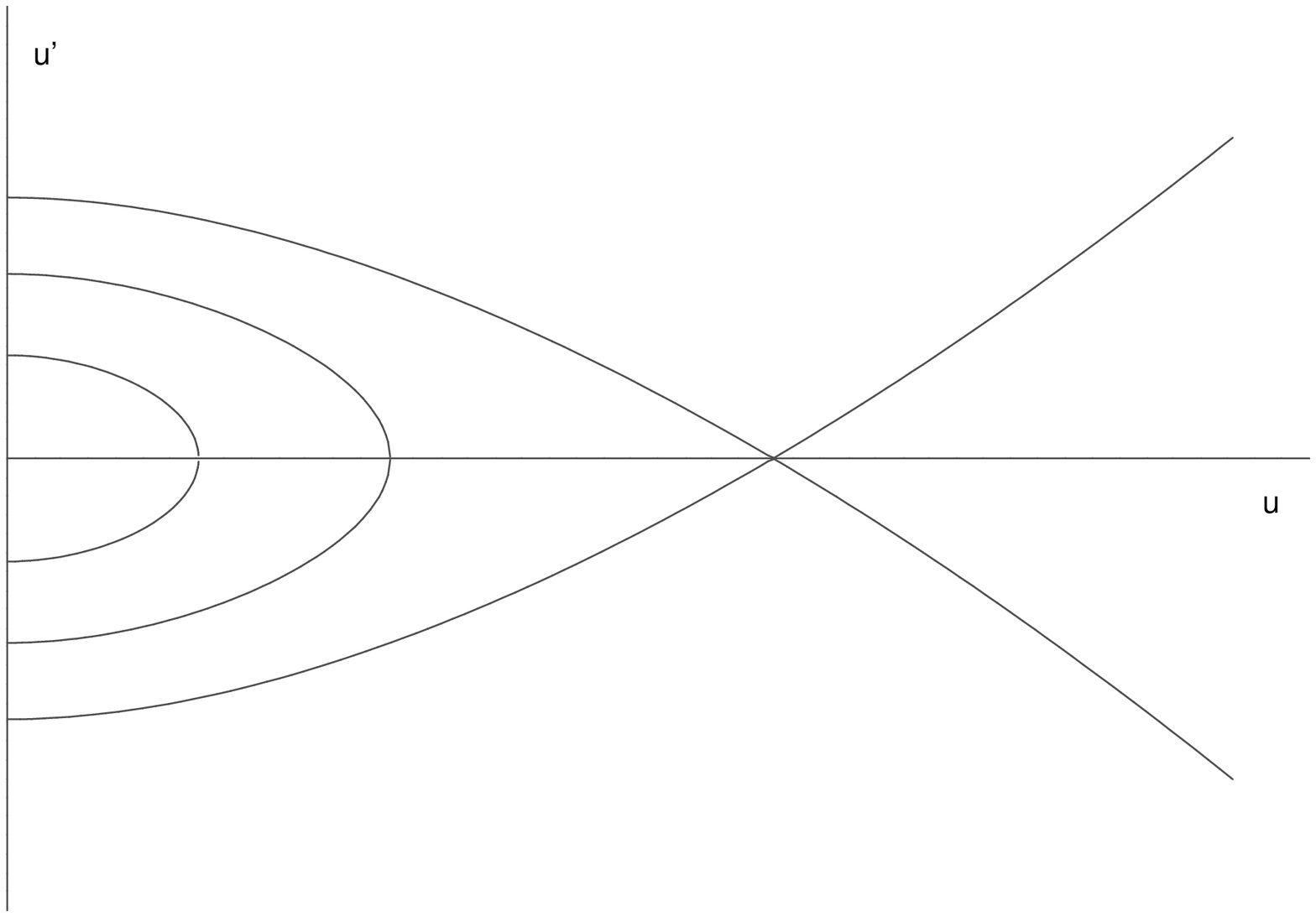}

One can see that there is no nonnegative solutions for any diffusion coefficient $D$ and any length $L$ of the environment.
On the other hand, if one use the Dirichlet boundary condition $u(0)=0=u(L),$ then a priori such solutions can exist. They will
correspond trajectories cutting axe $u'$ twice. One can compute the time (we use standard notions from the theory of dynamical systems
however, here, the independent variable is interpreted as a spatial one)  between two passes of this axe. Denote (after \cite{K})
by $F$ the real function $F(u)=u^2/2-u^3/3$ and let $\mu $ be the value of $u$ where the trajectory cut axe $u.$ Then this time
equals
$$L(\mu )=\sqrt{2D}\int _0^{\mu} \frac{du}{\sqrt{F(\mu )-F(u)}} =\int _0^1 \frac{\mu dz}{\sqrt{F(\mu )-F(\mu z)}}.$$
It is obvious that $\mu $ changes between $0$ and $1,$ that $L$ is
an increasing function of $\mu$ tending to $\infty $ as $\mu\to
1^{-} .$ One can also compute the limit
$$\lim _{\mu\to 0^+}L(\mu )=2\sqrt{D} \int _0^1 \frac{dz}{\sqrt{1-z^2}}=\pi \sqrt{D}.$$
It can be interpreted that there is exactly one nonconstant
steady-state solution iff the length of the environment $L$ is
greater than this limit $L_{KISS}=\pi\sqrt{D} $ called the KISS
size. If $L< \pi\sqrt{D} ,$ then there are only constant
steady-state  solutions $u\equiv 0$ and $u\equiv 1.$ From another
point of view, if the length $L$ is fixed, then nonconstant
time-independent solutions exist when the diffusion coefficient $D$
is sufficiently small.

Now, we study the case, when $\Om\subset \Bbb{R}^m$ with $m>1.$ The most interesting dimension from the biological point of view
is $m=2.$  $\Om =B(0,R)$ (the disk centered at $0$ with radius $R)$ is the simplest set and we can easily look for radial solutions
of (\ref{ell}). Our analysis can be repeated in larger dimension without troubles. If we denote by $u'$ the derivative with respect
to the radial coordinate, then we get the following boundary value problem:
\begin{equation}
u''+\frac{1}{r}u'=-\frac{1}{D}u(1-u)\quad {\rm in }\quad (0,R), \quad u'(0)=0=u'(R).  \label{steady2}
\end{equation}
For  Dirichlet's boundary condition $u|\p\Om =0, $  we have for radial solutions: $u'(0)=0=u(R).$ For both problems, we look
at solutions of the second order  ODE with initial values $u(0)=c,$ $u'(0)=0,$ where $c>0.$ It is easy to see that for $c>1$
the right hand side of ODE (\ref{steady2}) is positive, thus the solution increase to infinity in finite time and cannot satisfy
the Dirichlet nor Neumann conditions at $R.$ If $c\in (0,1),$ then the solution is concave near $0,$ $u$ decreases but, by comparison
with equation (\ref{steady}), slowler then in the last system. It means that, for any $c,$  $u(r_1)=0$ when $r_1>L_{KISS}/2.$
Thus, for Dirichlet's boundary condition, we have the same phenomenon as in dimension $1,$ but with larger KISS size.
Similarly, there is no nonconstant solutions for Neumann's problem. Below, we present numerical solutions for initial problems
$$u''+\frac{1}{r}u'=-\frac{1}{D}u(1-u),\quad u(0)=c,\quad u'(0)=0 $$
with several $c>0$ and $D=0.1.$

\includegraphics[angle=-90,width=12cm]{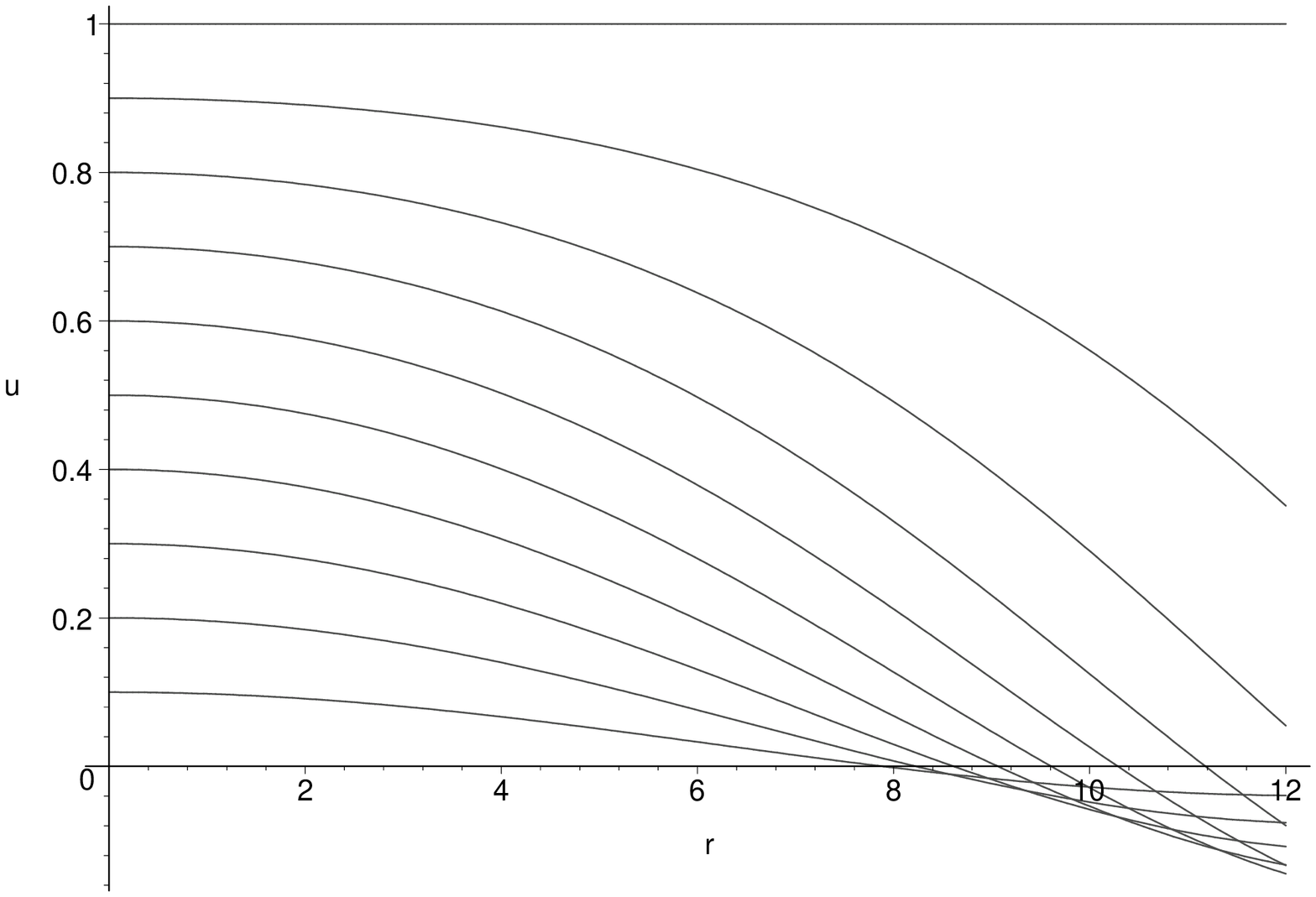}

\section{Steady-state solutions, two species}
\zero

For two species the situation is completely different. Consider the system
\begin{equation}           \label{2spec}
\left\{
\begin{array}{lll}
u_t=d_1\Delta u+u(1-u-bv),\;\;\; & \frac{\partial u}{\partial \nu}\left | _{\partial\Om}=0,\right. \;\;\; & u(\cdot ,0)=\varphi\\
v_t=d_2\Delta v+v(1-cu-v),\;\;\; & \frac{\partial v}{\partial \nu}\left | _{\partial\Om}=0,\right. \;\;\; & v(\cdot ,0)=\psi ,
\end{array}
\right.
\end{equation}
where interspecies-competition coefficients $b,c$ are positive constants, diffusion coefficients $d_1,d_2$ are positive and
$\Om $ is an open bounded set in $\Bbb{R}^m$ with smooth boundary. We study only nonnegative solutions. There are always three
equilibria constant in space and time:
$$P_0=(0,0),\quad P_u=(1,0),\quad P_v=(0,1),$$
and if $b,c<1$ or $b,c>1$, the fourth equilibrium
$$P_1=\left( \frac{1-b}{1-bc},\frac{1-c}{1-bc}\right) .$$
If $d_1=d_2=0,$ there is no diffusion and we have the standard
Lotka-Volterra ODE with the above equilibria. The dynamical system
given by this ODE is easily investigated: $P_0$ is its repeller,
$P_u$ is an attractor  if $c>1,$ $P_v$ is an attractor if $b>1,$
$P_1$ is global attractor (for nonnegative solutions) if $b,c<1.$
More exactly, if both $b$ and $c$ are greater than $1,$ then $P_1$
is a saddle point and its stable manifold $W^s(P_1)$ is the sum of a
heteroclinic trajectory from $P_0$ to $P_1$ and a trajectory from
infinity to $P_1.$ This manifold cuts the set $u,v\ge 0$ in two
sets: all trajectories starting from the set containing $P_u$ tend
to $P_u,$ trajectories starting from the second set tend to $P_v.$
The method of monotone dynamical systems enables us to state a
similar result for full parabolic system (\ref{2spec}):

(i) if $b,c<1,$ then  all solutions with nontrivial $\varphi , \psi  \ge 0$ tend to $P_1$ as $t\to \infty ;$

(ii) if $b<1$ and $c>1,$ then   all solutions with nontrivial $\varphi , \psi  \ge 0$ tend to $P_u$ as $t\to \infty ;$

(iii) if $b>1$ and $c<1,$ then   all solutions with nontrivial $\varphi , \psi  \ge 0$ tend to $P_v$ as $t\to \infty ;$

(iv) if $b,c>1,$ then the stable manifold $W^s(P_1)$ is again
one-dimensional, $P_u$ attracts solutions with $\varphi \ge u_0,$
$\psi\le v_0$ for some $(u_0,v_0)\in W^s(P_1),$ and $P_v$ attracts
solutions with  $\varphi \le u_0,$ $\psi\ge v_0$ for some
$(u_0,v_0)\in W^s(P_1).$

The proof can be found in \cite{S2}.

The last case is the most interesting since we do not know the dynamics for any $(\varphi ,\psi )$ not comparable in the above sense
with any function from $W^s(P_1).$ For example, some nonconstant equilibria can exist  and even they can be asymptotically stable.
However, we can use the following result due to Conway, Hoff and Sm\"oller \cite{CHS} in many cases:

Let $u$ be a solution to
$$u_t=D\Delta u+f(u),\quad  \frac{\partial u}{\partial \nu}\left | _{\partial\Om}=0,\right.  \quad u(\cdot ,0)=u_0,$$
where $u=(u_1,\ldots ,u_m),$ $D$ is a symmetric positively definite
matrix, $\Om $ is an open bounded subset of $\Bbb{R}^n$ with smooth
boundary, $f$ is a $C^2$ function and $u_0\in L^2(\Om ).$ Assume
that there exists a bounded positively invariant set $\Sigma\subset
\Bbb{R}^m$ such that $u_0$ takes values in $\Sigma .$ Denote by
$\lambda $ the first eigenvalue of $-\Delta $ with $\frac{\partial
u}{\partial \nu}\left | _{\partial\Om}=0\right .$ (here, one can
take $m=1,$ and we mean the first positive eigenvalue -- $0$ is also
eigenvalue), by $d$ the lowest eigenvalue of matrix $D,$ $M:=\sup
_{u\in\Sigma }||f'(u)||, $ and at last
$$\sigma := \lambda d-M.$$
Then there are four positive constants $c_i,$ $i=1,2,3,4$ such that
\begin{equation}       \label{1}
||\nabla u(\cdot , t)||_{L^2(\Om )}\le c_1{\rm e}^{-\sigma t},
\end{equation}
\begin{equation}       \label{2}
||u(\cdot , t)-\overline {u}(t)|| _{L^2(\Om )}\le c_2{\rm e}^{-\sigma t},
\end{equation}
where $\overline {u} $ is a spatial average of u and it satisfies $\overline {u}'=f(\overline {u})+g(t),$ $\overline {u}(0)=
\int _{\Om} u_0/\mu (\Om ),$ $g$ a is a function satisfying $|g(t)|\le c_3{\rm e}^{-\sigma t} ,$
\begin{equation}       \label{3}
||u(\cdot , t)- \overline {u}(t)|| _{L^{\infty}(\Om )}\le c_4{\rm e}^{-\sigma t}
\end{equation}
if $D$ is a diagonal matrix.

In our investigations, $\Sigma =\{(u,v):\; u,v\ge 0,\; u \le 1-bv\;
or \;  v\le 1-cu\}$ if $(\varphi ,\psi )$ takes values in this set
or $\Sigma$ is the smallest triangle with vortices $P_0,$ $(a,0),$
$(0,a)$ containing $(\sup \varphi , \sup \psi ).$ Since the first
set attracts all solutions, it contains all equilibria. If one has
$\sigma >0,$ then $\lim _{t\to\infty} \overline {u}(t) ={\rm const}$
and we the solution $u(x,t)$ tends to this constant as $t\to\infty
.$ It follows that there is no  nonconstant in space steady-state
solution. After easy though laborious computations one can find the
constant $M.$ We have
$$f'(u,v)=\left[
\begin{array}{ll}
1-2u-bv \;\;  & -bu\\
-cv & 1-2v-cu
\end{array}\right] ,$$
$$||f'(u,v)||^2=(b^2+c^2+4)(u^2+v^2)+4(b+c)uv-4(u+v)-2(cu+bv)+2$$
and the maximum of the last function on the set $\Sigma $ equals $2$ (it is reached at the origin). Thus, we have obtained

\begin{thm}
If both diffusion constants $d_1,d_2$ are sufficiently large, namely
$$\min (d_1,d_2)>\frac{\sqrt{2}}{\lambda}$$
where $\lambda $ is the first positive eigenvalue for $-\Delta $ with Neumann homogeneous condition, then there is no nonconstant
steady-state solution. For the case $\Om =(0,L)\subset \Bbb{R}$ as in the previous section, we have $\lambda =\frac{\pi ^2}{L^2}$
and we have no nonconstant equilibrium if
$$L<2^{-1/4}\pi \min (d_1,d_2).$$
\end{thm}
Compare this number with the KISS size from the previous section and notice that, here, a priori we have nonconstant equilibria for the Neumann
boundary condition. There is a numerically studied example of Matano and Mimura \cite{MM} where $\Om $ is a set in the plane
consisting two squares joined by a thin strip and $b,c>1$ such that there is a nonconstant positive equilibrium. Moreover, it is
asymptotically stable. If  domain $\Om $ is convex all equilibria are stable (see \cite{KW}), hence it is not surprising that
in this example the domain is such. On the other hand, one can understand ecological sense of the shape of the domain: in the first square
the first species wins, in the other the second one; the thin strip makes migrations between squares more difficult hence both species
can coexist.

Most results from the theory of monotone dynamical systems one used concern systems in ordered Banach spaces called SOP in the monograph
by H. Smith  \cite{S1}. Competitive parabolic systems are SOP if we use order:
$$(u,v)\prec (\check{u},\check{v})\Leftrightarrow u(x)\le \check{u}(x),\quad v(x)\ge \check{v}(x)\quad {\rm for\;\; any}\;\; x.$$
The semiflow generated by (\ref{2spec}) preserves this order, i.e.
if $(\varphi ,\psi)\prec (\check{\varphi},\check{\psi}),$ then
$(u(\cdot ,t), v(\cdot ,t))\prec (\check{u}(\cdot
,t),\check{v}(\cdot ,t))$ for any $t>0.$ It enables us to get some
information on $\omega $-limit sets of our system. This choice of
the order is typical for competition of two species and can be
explain qualitatively by ecological arguments: if the first species
dominates the second one and one increases the  population of the
first species and decreases of the second one then the domination
conserves.

\section{Steady-state solutions, many species}
\zero

For more than two species the situation is much more complicated. If
$n=3$ then, roughly speaking the first species can dominate the
second one, the second one can dominate the third one and this last
species can dominate the first one. The simplest mathematical model
of this case is given by May and Leonard \cite{ML} for ODE:
\begin{equation} \label{MayL}
\left\{
\begin{array}{l}
\dot{x}= x(1-x-\alpha y-\beta z)\\
\dot{y}= y(1-\beta x-y-\alpha z)\\
\dot{z}= z(1-\alpha x-\beta y-z)
\end{array}
\right.
\end{equation}
 with $\alpha ,\beta >0,$ $\alpha + \beta   =2.$ The typical $\omega$-limit set is a limit cycle but there are three heteroclinic
 trajectories joining three equilibria $(1,0,0),$ $(0,1,0)$ and $(0,0,1).$ We have studied a slightly more general case in \cite{JP}
 however the behavior of the system is very similar.

 Consider the parabolic system for three species:
 \begin{equation} \label{sys3}
 \begin{array}{l}
 U_t=D\Delta U+ f(U), \\
 \frac{\p U}{\p\nu}\left |_{\p\Om}=0,\right.\\
 U(\cdot ,0)=\varphi   ,
 \end{array}
 \end{equation}
 where $U=(u,v,w),$  $D={\rm diag} (d_1,d_2,d_3),$
 \begin{equation} \label{nonl}
 f(U)=\left [
 \begin{array}{l}
  u(1-a_1u-b_1v-c_1w)\\
  v(1-a_2u-b_2v-c_2w)\\
  w(1-a_3u-b_3v-c_3w)
 \end{array}
 \right ],
\end{equation}
all constants in the above formulas $d_i,$ $a_i,$ $b_i,$ $c_i$ are positive, $x\in\Om\subset \Bbb{R}^m.$

We study the system (\ref{sys3}) and its spatially homogeneous ODE system
\begin{equation} \label{ODE}
U'=f(U)
\end{equation}
in the open set
$${\cal D}:=\{ (u,v,w):\; u,v,w> 0\}$$
which is obviously invariant and the same is true for its closure $\overline{\cal D}.$ Fixed points of the system (\ref{ODE}) can be easily found
-- four of them always lie in
$\overline{\cal D}:$
$$P_0=(0,0,0),\quad P_u=(a_1^{-1},0,0), \quad P_v=(0,b_2^{-1},0) , \quad P_w=(0,0,c_3^{-1}),$$
next four are
$$P_{uv}=\left(\frac{b_2-b_1}{a_1b_2-a_2b_1},\frac{a_1-a_2}{a_1b_2-a_2b_1},0\right),\quad
P_{uw}=\left(\frac{c_3-c_1}{a_1c_3-a_3c_1},0,\frac{a_1-a_3}{a_1c_3-a_3c_1}\right),$$
$$\qquad P_{vw}=\left(0,\frac{c_3-c_2}{b_2c_3-b_3c_2},\frac{b_2-b_3}{b_2c_3-b_3c_2}\right)$$
and $P_1=(\alpha ,\beta ,\gamma )$ being the solution of the equation
$$M\left[ \begin{array}{l}   u\\v\\w \end{array}\right]=  \left[ \begin{array}{l}   1\\1\\1 \end{array}\right]
\quad {\rm where} \quad M= \left[ \begin{array}{lll}   a_1 & b_1 & c_1\\ a_2 & b_2 & c_2\\ a_3 & b_3 & c_3  \end{array}\right].$$
We assume that $W:=\det M \neq 0$ what means that none of two of planes
$$H_i:\quad a_iu+b_iv+c_iw=1,\quad i=1,2,3$$
are parallel. It is easy to see that
$f'(P_0)=I $
and  its unique eigenvalue $1$ is positive, thus $P_0$ is a source for (\ref{ODE}).
Similarly
$$f'(P_u)=\left[ \begin{array}{ccc}  - 1 & -b_1 /{a_1} &  -c_1 /a_1 \\ 0 & 1-a_2/a_1 & 0 \\ 0 & 0 & 1-a_3/a_1  \end{array}\right]  $$
and if $\max (a_2,a_3)< a_1,$ then one eigenvalue $-1$ is negative and two remaining ones are positive. The corresponding eigenspaces
are: $\{(u,0,0):\; u\in\Bbb{R}\}$ and  $\{(0,v,w):\; v,w\in\Bbb{R}\}.$ The stable manifold of $P_u$ is $1$-dimensional:
$$\{(u,0,0):\; u>0\}$$
and by the Stable and Unstable Manifold Theorem all trajectories of (\ref{ODE}) except starting in the stable manifold cannot tend to $P_u$ as $t\to +\infty.$
Similar arguments work for $P_v$ and $P_w$ under conditions:
$$\max (b_1,b_3)<b_2,\quad \max (c_1,c_2)<c_3.$$
 In \cite{JP}, we proved even more:
\begin{lem} If
\begin{equation} \label{ineq}
\min(a_2,a_3)<a_1,\quad \min(b_1,b_3)<b_2,\quad \min(c_1,c_2)<c_3,
\end{equation}
then $P_u,$ $P_v$ and $P_w$ do not belong to the $\omega$-limit set $\omega (P)$ of any point $P\in {\cal D}.$
\end{lem}
Dividing the set $\overline{\cal D} $ into three pieces:
$${\cal D}_+:=\{(u,v,w)\in \overline{\cal D} :\; \min _i (a_iu+b_iv+c_iw)>1\},$$
$$A:=  \{(u,v,w)\in \overline{\cal D} :\; \min _i (a_iu+b_iv+c_iw)\le 1\le\max _i(a_iu+b_iv+c_iw)\},$$
$${\cal D}_-:=\{(u,v,w)\in \overline{\cal D} :\; \max _i (a_iu+b_iv+c_iw)<1\}$$
(${\cal D}_+$ (resp. ${\cal D}_-$) is the set of points sitting under (resp. over) all three planes $H_i,$ $i=1,2,3,$ $A= \overline{\cal D}\setminus ({\cal D}_+\cup {\cal D}_-)$)
we got \cite{JP} following result for (\ref{ODE}):
\begin{lem}
The set $A$ is positively invariant and all trajectories in ${\cal D}$ eventually come into $A.$ Moreover, $A$ contains any compact invariant set that contains no fixed points.
\end{lem}
Thus, $\omega (P)\subset A$ for any $P\in {\cal D}.$ Notice that $P_u,$ $P_v,$ $P_w\in A$ and similarly $P_{uv},$ $P_{uw},$ $P_{vw}$ if they belong to $\overline{\cal D}.$

The last fixed point $P_1=(\alpha ,\beta ,\gamma )$ can lie in ${\cal D}$ (and then in $A$) or outside this cone. The following theorem excludes the existence of a periodic trajectory for (\ref{ODE}) in ${\cal D}$ if
$P_1\in {\cal D}.$
\begin{lem} (\cite{S1}, p. 44, Prop. 4.3)
Let $\Gamma $ be a nontrivial periodic trajectory of a competitive system in ${\cal D}\subset \Bbb{R}^3$ and
$$\Gamma \subset [p,q]:=\{\xi: \; p_i\le\xi _i\le q_i,\; i=1,2,3\}\subset D.$$
Then the set $K$ of all points $x$ which are not related to any point $y\in \Gamma $ (relation $x\le y$ means $x_i\le y_i$ for any $i$)
has two components, one of them is bounded and contains a fixed point.
\end{lem}
Hence, if we want to have a nontrivial periodic trajectory, then $P_1$ must belong to ${\cal D}$ and we have two options:

(i) $W>0$ and three other determinants
$$W_u:=\det \left[ \begin{array}{lll}   1 & b_1 & c_1\\ 1 & b_2 & c_2\\ 1 & b_3 & c_3  \end{array}\right],$$
$$W_v:=\det \left[ \begin{array}{lll}   a_1 & 1 & c_1\\ a_2 & 1 & c_2\\ a_3 & 1 & c_3  \end{array}\right],$$
$$W_w:=\det \left[ \begin{array}{lll}   a_1 & b_1 & 1\\ a_2 & b_2 & 1\\ a_3 & b_3 & 1  \end{array}\right]$$
are positive or

(ii) $W<0$ and the above three determinants are negative.

The main result of \cite{JP} is the following
\begin{thm}        \label{main}
Assume (\ref{ineq}). Let all four determinants $W,$ $W_u,$ $W_v,$ $W_w$ be positive and $p:=a_1\alpha +b_2\beta +c_3\gamma -1<0.$
Then, for any point $P\in {\cal D}$ that does not belong to the half-line
$$\left\{
\begin{array}{ll}
u=\alpha s, &  \\
v=\beta s,\quad & s>0,\\
w=\gamma s &
\end{array}
\right .$$
  the $\omega$-limit set $\omega (P)$ for (\ref{ODE}) is
a periodic trajectory. For $P$ from this half-line,   $\omega (P)=P_1.$
\end{thm}
If we combine this theorem with the result \cite{CHS} cited in the previous section and a classical result of L. Markus \cite{Ma}
that $\omega$-limit sets of (\ref{sys3}) are the same as corresponding (\ref{ODE}), we obtain
\begin{thm}
Suppose that inequalities (\ref{ineq})  hold, four determinants $W,$ $W_u,$ $W_v,$ $W_w$ are  positive,
$p<0$ and all diffusion coefficients $d_i$ are sufficiently large, then for any continuous, nonnegative function
$\varphi$ with all coordinates  nontrivial, either
$$\lim _{t\to\infty }U(x,t)=P_1$$
or  there exists a periodic function $\tilde{U}:\Bbb{R}\to A\subset {\cal D}$ such that
$$\lim _{t\to\infty }|U(x,t)-\tilde{U}(t)|=0$$
and in both cases limits are uniform with respect to $x.$  The first case takes place only if
$\int _{\Om}\varphi $ belongs to the half-line starting from the origin.
\end{thm}
The condition for the diffusion coefficients is given by:
$$\lambda _1\min \{ d_i:\; i=1,2,3\}>\sup \{||f'(U)||:\; U\in A\}$$
where $\lambda _1$ is the first positive eigenvalue of $-\Delta $ with Neumann's boundary condition, but numerical experiments
show that the assertion of the above theorem holds also for essentially larger these coefficients and even if only one of them
is great and two others are very small. However, for all coefficient being small there exists spatially heterogeneous steady-state
solutions which attracts other solutions. Obviously, these functions take values in positively invariant subset $A.$

If the number $p$ defined in Theorem \ref{main} is positive, then all eigenvalues of $f'(P_1)$ have negative real parts --
calculations in \cite{JP} show that the characteristic polynomial equals
$$Q(\lambda)=(\lambda +1)( \lambda ^2+p\lambda  +\alpha\beta\gamma W).$$
It follows that
$$\lim _{t\to\infty }U(x,t)=P_1$$
 uniformly in $x\in\Om $ for any starting point $\varphi $ as in the last theorem.  The case $p=0$ is difficult for investigations
 since the behavior of the system near $P_1$ cannot be found by the linearization; $P_1$ is not hyperbolic. In the special case
 from the paper of May and Leonard \cite{ML} $a_1=b_2=c_3=1,$ $b_1=c_2=a_3,$ $c_1=a_2=b_3$ if $p=0$ there are two first integrals
 of (\ref{ODE}) and trajectories can be found apparently. The whole triangle spanned by $P_u,$ $P_v$ and $P_w$ is fulfilled by
 limit cycles around $P_1.$ For sufficiently small diffusion coefficients these limit cycles describe the asymptotic behavior of
 the system (\ref{sys3}).

 If $P_1$ does not belong to $\cal D,$ then all solutions  tends to one of the others steady-state solutions at least for small
 diffusion coefficients, when there is no nonconstant steady-state solutions due to the above mentioned arguments. It means that
 at least one of the species extincts.

 \section{Some numerical simulations}
 \zero

 The most interesting case considered in the previous section is presented in the theorem: almost all solutions of (\ref{sys3}) tend
 to periodic functions given by the system of ODEs (\ref{ODE}). We can investigate numerically equation
 (\ref{sys3}) with $\Om =(0,1)\subset \Bbb{R}$ and $f$ of the form from  the previous section. Put
 matrix $M$ of competition coefficients
 $$ M= \left[ \begin{array}{lll}   2 & 1.1 & 3.1\\ 3.1 & 2 & 0.9\\ 0.95 & 2.9 & 2  \end{array}\right].$$
 This choice ensure the assumptions of Theorem 3 hold. The first positive eigenvalue of this degenerate Laplacian equals
 $\lambda _1= \pi ^2 $ and one can compute the maximal value of the norm of the derivative
 $||f'(U)||$ on the set $A:$
 $$||f'(U)||\le \sqrt{3}.$$
 Thus the critical value of diffusion coefficients is  $\tilde{d}:=\frac{\sqrt{3}}{\pi ^2}.$ For $\min \{ d_i:\; i=1,2,3\}>
 \tilde{d},$ almost all trajectories tend to periodic, spatially constant functions. Since the theorem of Conway, Hoff and Sm\"oller
 gives only the sufficient condition, it is not surprising that even for smaller $d_i$'s the assertion is true.

 We have found numerically (by using Maple 10) solutions of (\ref{sys3}) with the above matrix $M$ and
 $$d_1=10^{-3},\quad d_2=2\cdot 10^{-3},\quad d_3=0.5\cdot 10^{-3}$$
 $$\varphi (x)=\left[ 6x^2(1-x)^3,x^4(1-x)^2,2x^3(1-x)^2\right].$$
 The type of this initial function is natural if we are seeking classical solutions: the normal derivative
 of the initial function should vanish
 at both boundary points $x=0,$ $x=1.$  Below, we present the plot of the graph of the second coordinate of $U$ for three
 values of $x:$ $0.1,$ $0.5$  and $0.9$ as the function of time $t$ for the range $[0,100]:$

 \includegraphics[angle=-90,width=12cm]{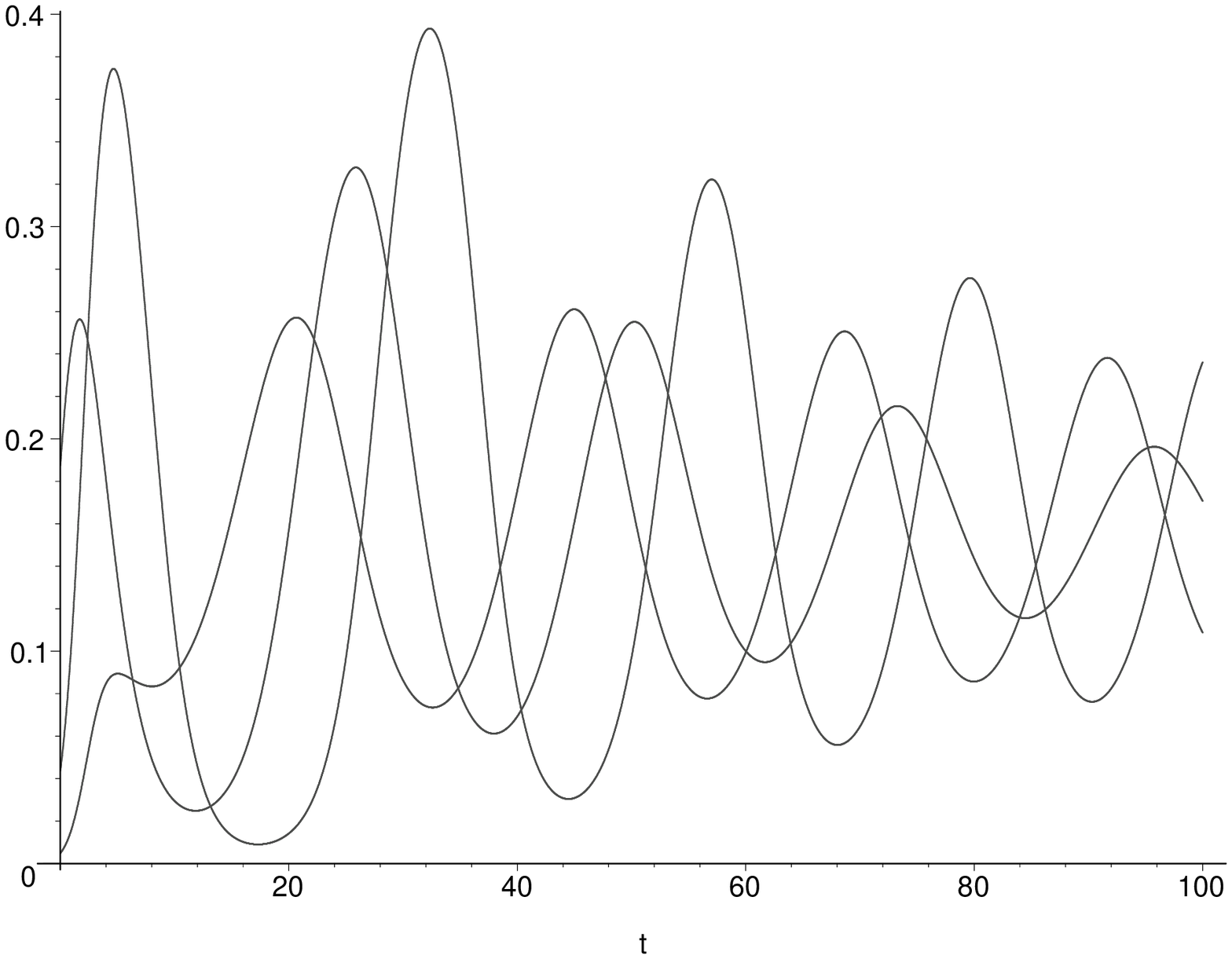}

The plots for three values of $x$ are  completely different, hence the solutions of the parabolic system
tend as $t\to\infty$ to functions which are not spatially constant. Nevertheless, they seem to be periodic as functions of time. Compare this plot
with the plot of the second coordinate $v$ of the solution to ODE (\ref{ODE}) with initial point $(0.1,0.0095238,0.0333333)$
which is the spatial average of $\varphi .$

\includegraphics[angle=-90,width=12cm]{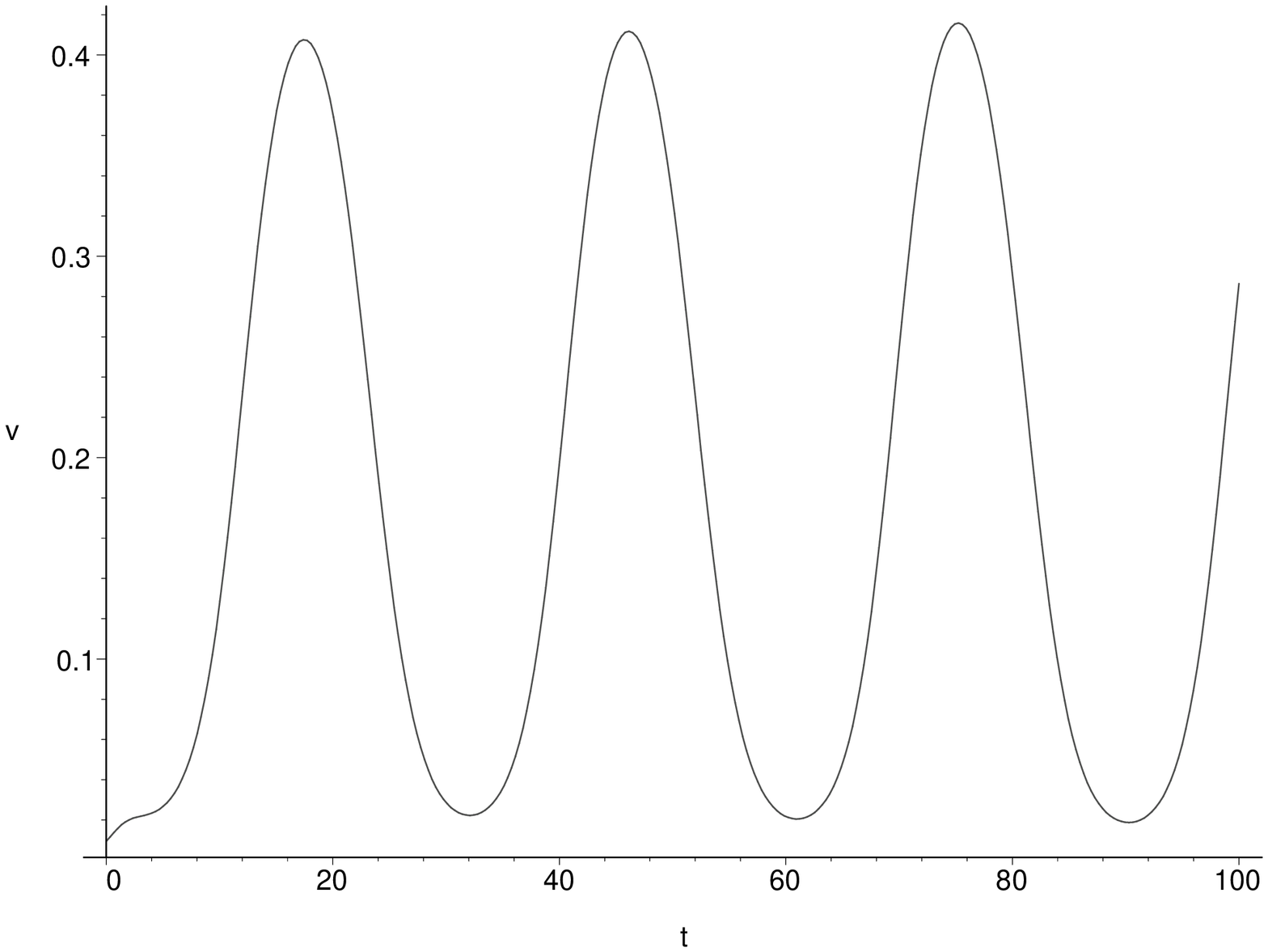}

\section{Stability of steady-state solutions}
\zero The stability analysis in the linear approximation for
constant solutions of the parabolic system is standard (see, for
example \cite{Sm,I}). Consider (\ref{sys3}) not necessarily with $U$
taking values in $\Bbb{R}^3$ but in $\Bbb{R}^n.$ If $P$ is a zero of
$f,$ then $U\equiv P$ is a solution both (\ref{sys3}) and
(\ref{ODE}) and it is asymptotically stable for both systems again
if  all eigenvalues of $f'(P)$ have negative real parts. Below, we
shall study the analogous problem for periodic solutions of
(\ref{ODE}).

Let $t\mapsto U(t)$ be such a solution. Denote by $A(t)=f'(U(t)),$ by $\lambda _k,$ $k=0,1,2,\ldots ,$ the sequence of all
eigenvalues of $-\Delta $ with Neumann boundary conditions ($\lambda _0=0<\lambda _1\le \lambda _2\le \ldots )$ and by $x\mapsto e_k(x)$
corresponding eigenfunctions. Since the spectrum is discrete and the operator is self-adjoint, the eigenfunctions form a complete
orthonormal system in $L^2(\Om ).$ One can solve the linearized problem
\begin{equation}
U_t=D\Delta U+A(t)U,\quad \frac{\p U}{\p\nu}\left |_{\p\Om}=0,\right. \quad U(\cdot ,0)=\varphi       \label{lin}
\end{equation}
by using the Fourier method.  If we put
$$U(x,t)=\sum _{k=0}^{\infty } e_k(x) g_k(t)  $$
in  (\ref{lin}), we have Neumann's condition satisfied and vector-valued function $g_k$ should be a solution of the initial
problem:
\begin{equation}      \label{ini}
   g_k'(t)=\left ( -{\lambda _k} D+A(t)\right )g_k(t),\qquad g_k(0)=c_k
\end{equation}
  for any $k,$ where $c_k$ is a coefficient of the Fourier expansion of $\varphi $ with respect to the orthonormal system $\{ e_k:\; k=0,1,\ldots \},$
  i.e.
  $$c_k=\int _{\Om }\varphi e_k.$$
Since matrix $D$ is diagonal, it commutes with $A(t)$ and
multipliers (Floquet theory) of $-{\lambda _k} D+A(\cdot)$ are of
the form $\exp(-\lambda _k d)\varrho ,$ where  $d$ belongs to the
interval $[\min \{d_i\},\max \{d_i\}]$ and $\varrho $ is a
multiplier for matrix $A(\cdot ).$ Hence $g_k$ decays exponentially
as $t$ tends to $+\infty ,$  if $|\varrho |< \exp(\lambda _k \min
\{d_i\} )$  for any multiplier $\varrho .$ It is well known that one
of multipliers for $A$ equals 1, thus this inequality cannot hold
for $k=0$ as $\lambda _0=0.$ But, in spite of this, if the remaining
multipliers of $A$ sit in the open unit disc and $1$ is simple, then
solution $U$ of (\ref{sys3}) is orbitally  asymptotically stable --
see \cite{Sm} Theorem 8.2.3, p.251, i.e. there exists a neighborhood
of the periodic orbit $\Gamma :=\{ U(t):\; t\in \Bbb{R}\}$ such that
solutions $U_1$ with initial function $\varphi $ taking values in
this nhbd tend to $U$ in the sense
$$\lim _{t\to +\infty}||U_1(\cdot ,t)-U(t+t_0)||=0,$$
where the above norm means usual one in $H^1$ and $t_0$ is a number depending on $\varphi $ called asymptotic phase. For $1$-dimensional
domains $\Om ,$ one has the canonical embedding of $H^1$ and the above limit is uniform in $x.$ For more natural environments ($\Om\subset \Bbb{R}^2),$
we cannot use this argument to get the uniform limit.

\section{Concluding remarks}
\zero
Recently, a lot of important papers on spatial heterogeneity models of competiting species \cite{Du1,Du2,HLM1,HLM2,LGM,L}, for instance.
They consider mutual interplay of diffusion and competition which gives many interesting phenomena for such systems.
However, most of results are obtained for two species but spatial heterogeneity is included in the model -- some coefficients
depend on variable $x.$ It seems that some new effects  can be obtained if we study interaction of more than two species.
The present author hope that some possible directions of such investigations have been indicated above.

\end{document}